\documentclass[a4paper,11pt]{elsarticle}
\usepackage{amsfonts}
\usepackage{amssymb}
\usepackage{amsmath}
\usepackage{amsthm}

\usepackage[utf8]{inputenc}
\usepackage[french]{babel}
\usepackage[T1]{fontenc}

\usepackage{geometry}
\theoremstyle{definition}
\newtheorem{definition}{Definition}

\begin{document}

\begin{frontmatter}

  \title{On near orthogonality of the Banach frames of the wave packet spaces}
  
  \author{Dimitri Bytchenkoff $^{1, 2, *}$}
  
  \address{$^{1}$Österreichische Akademie der Wissenschaften, Institut für Schallforschung, Wohllebengasse 12-14, 1040 Wien, Österreich}
  \address{$^{2}$Laboratoire d'Energ{\'e}tique et de M{\'e}canique Th{\'e}orique et Appliqu{\'e}e, Universit{\'e} de Lorraine, 2 avenue de la For{\^e}t de Haye, 54505 Vandoeuvre-l{\`e}s-Nancy, France}
   
 \begin{abstract}
In solving scientific, engineering or pure mathematical problems one is often faced with a need to approximate the function of a given class by the linear combination of a preferably small number of functions that are localised one way or another both in the time and frequency domain. Over the last seventy years or so a range of systems of thus localised functions have been developed to allow the decomposition and synthesis of functions of various classes. The most prominent examples of such systems are Gabor functions, wavelets, ridgelets, curvelets, shearlets and wave atoms. We recently introduced a family of quasi-Banach spaces -- which we called $wave$ $packet$ $spaces$ -- that encompasses all those classes of functions whose elements have sparse expansions in one of the above-mentioned systems, supplied them with Banach frames and provided their atomic decompositions. Herein we prove that the Banach frames and sets of atoms of the wave packet spaces -- which we call $wave$ $packet$ $systems$ -- are indeed well localised or, more specifically, that they are near orthogonal. We believe that good localisation of the wave packet systems -- a property of paramount importance that can't be taken for granted even for frames in Hilbert spaces, let alone Banach spaces -- will pave the way, among other things, for their use for representing Fourier integral operators on Banach spaces by sparse and well structured matrices by the Galerkin method. This, in its turn, should allow one to design efficient computer programmes for solving corresponding operator equations on two-dimensional manifolds. \\
\newline

\noindent $MSC:$ 42B35, 42C15, 42C40
   \end{abstract}
      
   
\end{frontmatter}

\let\thefootnote\relax\footnotetext{* Corresponding author. E-mail address: dimitri.bytchenkoff@univ-lorraine.fr (D. Bytchenkoff)}

\section{Introduction}

\noindent In science \cite{Rodts_2013, Rodts_2011, Reisenhofer_2016}, engineering \cite{Chambolle_1998, Luisier_2007, Kutyniok_2017, Adcock_2015} and mathematics \cite{Liang_2012, Voigtlaender_2017, Frazier_1989, Groechenig_2006, Meyer_1990, Kato_2019, Fell_2016, Kutyniok_2009, Candes_2005}, functions often need to be decomposed into or synthesised from those localised one way or another both in the time and frequency domain. For this purpose a number of systems of localised real functions of two real variables have been designed and successfully used over the last seventy years or so. The most important of them are Gabor functions \cite{Groechenig_2001}, wavelets \cite{Daubechies_1992, Triebel_2006, Frazier_1985, Frazier_1989, Bownik_2005, Borup_2007, Bytchenkoff_2020}, $\alpha$-modulation functions \cite{Groebner_1992}, ridgelets \cite{Candes_1999, Grohs_2012}, curvelets \cite{Candes_2002}, shearlets \cite{Kutyniok_2011, Kittipoom_2012} and wave atoms \cite{Demanet_2009}.

Two essential properties of each of these systems, except for the anisotropic wavelets, are determined by two parameters known as $\alpha$ and $\beta$. The parameter $\alpha$ defines how the length $(1 + \vert \xi \vert)^\alpha$, so to speak, of the Fourier transform of the element of the system along its radial axis of symmetry depends on the absolute value of the frequency $\xi$ at which it is localised. For instance, the length of the Gabor functions, whose $\alpha$ equals zero, is independent of the frequency, while that of wavelets, whose $\alpha$ equals one, is proportional to it. That being so, expanding a given function in a series of Gabor functions will result in identifying higher and lower frequencies present in the function with the same resolution,
while in its expansion in a series of wavelets lower frequency will be much better resolved than higher ones. The parameter $\beta$ determines the number $2^{(1-\beta)j}$ of elements of the systems whose Fourier transforms appear at a given frequency scale, defined by the index $j$, and differ from each other only in orientation. This, in its turn, defines directional resolution of different frequencies present in the function that can be achieved by expanding it into a series of a given type of localised functions. For instance, Gabor functions, for whose $\beta$ equals zero, have equally good direction resolution at all frequency scales, while wavelets, for whose $\beta$ equals one, do not allow one to distinguish different directions in the frequency domain at all.

The spaces of functions whose elements are characterised by their sparse decompositions into localised functions that belong to one or another of the systems mentioned above are worth studying in their own right. Not long ago we introduced a family of quasi-Banach spaces -- which we called $wave$ $packet$ $spaces$ -- that encompasses all these spaces, studied their properties, equipped them with Banach frames -- which we called $wave$ $packet$ $systems$ -- and provided their atomic decompositions \cite{Bytchenkoff_2019}. We adopted this term from \cite{Cordoba_1978} as we believe that our wave packets unify in themselves the quasi-Banach frames \cite{Groechenig_1991} mentioned above as the wave packets in \cite{Cordoba_1978, Hernandez_2002, Hernandez_2004, Labate_2004} unified frames for separable Hilbert spaces \cite{Duffin_1952}. 
Efficiency of the approximation of functions \cite{Grohs_2014} of various classes in terms of a combination of elements of an appropriately chosen wave packet system and their potential usefulness for discretisation of bounded linear operators \cite{Balazs_2017}, in general, and Fourier integral operators \cite{Cordoba_1978, Candes_2004}, in particular, by the Galerkin method can be expected to depend crucially on how well these elements are localised, not least because, the Fourier integral operators are known to essentially transform well localised wave packets into each other \cite{Cordoba_1978}. Knowing, say from \cite{Labate_2004}, that good localisation of frames even for Hilbert spaces should not be taken for granted, we give herein a precise and thorough answer to this question.

To make this manuscript self-contained, we first remind the definition of the wave packet spaces, in Section 2. We then slightly reformulate the definition of the the wave packet systems so that they will, on the one hand, contain no two identical wave packets and, on the other hand, be both Banach frames and sets of atoms for the wave packet spaces as long as their parameters $\alpha$ and $\beta$ satisfy the condition $0 \leqslant \beta \leqslant \alpha \leqslant 1$, discuss their structure, remind the notion of good localisation of the set of functions in terms of near orthogonality of its elements and prove it for the wave packet systems, in Section 3. Final conclusions are drawn in Section 4. The list of notations used throughout this manuscript is provided in Section 5.

\section{Structure of wave packet spaces}

\noindent The wave packet spaces are quasi-Banach spaces defined by the decomposition method \cite{Feichtinger_1985} and, as such, are made of three basic building blocks, namely an almost-structured covering $\mathcal{Q}$ of a Lebesgue-measurable subset of the frequency plane, a regular partition of unity and a $\mathcal{Q}$-moderate weight. Here are their definitions as well as those of all necessary axillary notions.

\begin{definition}
The set $\mathcal{Q} = \{ Q_i \}_{i \in I}$ is called an $almost$ $structured$ $covering$ of a Lebesgue-measurable subset $O$ of $\mathbb{R}^2$ if
\begin{enumerate}
  \item the number of elements in the sets $\left\{i' \in I : Q_{i'} \cap Q_i \neq \emptyset \right\}$ is uniformly bounded for all $i \in I$;
    
  \item there exists a set $\{ T_i \bullet + b_i \}_{i \in I}$ of invertible affine-linear maps and finite sets $\{ Q'_n \}_{n=1}^N$ and $\{P_n\}_{n=1}^N$ of non-empty open and bounded subsets $Q'_i$ and $P_i$ of $O$ such that
   \begin{enumerate}                 
      \item $\overline{P_n} \subset Q'_n$ for $1 \leqslant n \leqslant N$;
      \item for each $i \in I$ there exists such an $n_i \in \{ 1, ..., N \}$ that $Q_i = T_i \, Q'_{n_i} + b_i$;        
      \item there exists such a constant $C_1 > 0$ that $\|T_i^{-1}T_{i'}\| \leqslant C_1$ for all such
          $i$ and $i' \in I$ that $Q_i \cap Q_{i'} \neq \emptyset$; and
      \item $\bigcup_{i \in I} (T_i P_i + b_i) = O$.
    \end{enumerate} 
\end{enumerate}
\label{eq:Structuredness}
\end{definition}

\begin{definition}
First, let $\epsilon \in (0, 1/32)$ and
\begin{equation}
Q := (-\epsilon, 1+\epsilon) \times (-1-\epsilon, 1+\epsilon)
 \hspace{0.25cm};
\label{eq:Mutter-Kachel}
\end{equation} second, let $0 \leqslant \beta \leqslant \alpha \leqslant 1$, $N:=10$, $j^{\max} \in \mathbb{N}$ and
\begin{equation}
I_0 := \{ (0, 0, 0) \} \cup \{(j, m, l) \in \mathbb{N} \times \mathbb{N}_0 \times \mathbb{N}_0 :
 j \leqslant j^{\max}, \text{ } m \leqslant m_j^{\max} \text{ and } l \leqslant l_j^{\max} \}
\label{eq:Index-Menge}
\end{equation} where
\begin{equation}
 l_j^{\max}
    := \lceil
           N \cdot 2^{\left(1-\beta\right)j}
         \rceil -1
          \hspace{0.25cm},
\label{eq:l-Index-Grenzen}
\end{equation}
\begin{equation}
 m_j^{\max}
    := \begin{cases}
           \lceil 2^{\left(1-\alpha\right)j-1} \rceil - 1
          \, 
         & \text{if } \theta_{jl} = \theta_{j+1, l'} \text{ for an }  \mathbb{N}_0 \ni l' \leqslant l_{j+1}^{\max} \text{ or if } j=j^{\max} \, \\      
         \lceil 2^{\left(1-\alpha\right)j-1} \rceil \, 
         &
           \text{if } \theta_{jl} \neq \theta_{j+1, l'} \text{ for any }  \mathbb{N}_0 \ni l' \leqslant l_{j+1}^{\max} \,
       \end{cases}
\label{eq:m-Index-Grenzen}
\end{equation} and
\begin{equation}
\theta_{jl} := 2\pi \cdot \frac {{2}^{ \left( \beta -1 \right) j}} {N} \cdot l
\hspace{0.25cm};
\label{eq:theta_jl}
\end{equation}
finally, let $Q_0 := Q_{000} := B_4(0)$, where $B_4(0)$ stands for the Euclidean circle of radius four with its centre in the origin, and
\begin{equation}
Q_i := Q_{jml} := R_{jl} (A_j Q + B_{jm})
\label{eq:Tochter-Kacheln}
\end{equation} where
\begin{equation}
A_j = \left( \begin{array}{ cc }
2^{\alpha j} & 0 \\
0 & 2^{\beta j} \\
\end{array} \right)
\text{ , }
 \hspace{0.25cm}
B_{jm} = \left( \begin{array}{ cc }
2^{j-1} + m \cdot 2^{\alpha j} \\
0  \\
\end{array} \right)
\label{eq:A_j_und_B_jm}
\end{equation} and
\begin{equation}
R_{jl} := \left( \begin{array}{ cr }
\cos \theta_{jl} & - \sin \theta_{jl} \\
\sin \theta_{jl} & \cos \theta_{jl} \\
\end{array} \right)
\hspace{0.25cm}
\label{eq:R_jl}
\end{equation} for all $i = (j, m, l) \in I_0 \setminus \{ (0, 0, 0) \}$. The set $\mathcal{Q}^{(\alpha,\beta)} := \{ Q_i \}_{i \in I_0}$
is called $wave$ $packet$ $covering$ of the subset $O$ of $\mathbb{R}^2$, where $O$ is identical with $\mathbb{R}^2$ if $j^{\max}$ is infinite or such a polygon that $B_{R_1} (0) \subset O \subset B_{R_2} (0)$ where $R_1 := 2^{j^{\max}}$ and
\begin{equation*}
R_2 := \frac{2^{j^{\max}}}{ \cos \left( \pi \cdot \frac{2^{(\beta - 1) j^{\max}}}{N} \right)}
\end{equation*} if $j^{\max}$ is finite, the precise definition of $O$ being of little interest.
\label{eq:Ueberdeckung}
\end{definition}

The proof that $\mathcal{Q}^{(\alpha,\beta)}$ is indeed a covering of $O$ and that it is almost-structured is almost identical with that of Lemma 3.2 and 5.1 of \cite{Bytchenkoff_2019} respectively.

\begin{definition}
Let $\mathcal{Q} = \{ Q_i \}_{i \in I}$ and $\{ T_i \bullet + b_i \}_{i \in I}$  be an almost structured covering of a Lebesgue-measurable subset $O$ of $\mathbb{R}^2$ and the set of invertible affine-linear maps associated with it respectively. The set of functions $\Phi = \{ \phi_i \}_{i \in I}$ is called $regular$ $partition$ $of$ $unity$ $subordinate$ $to$ $\mathcal{Q}$,
  if
  \begin{enumerate}
    \item $\phi_i \in C_c^\infty (\mathbb{R}^2)$ with
          $\operatorname{supp} \phi_i \subset Q_i$ for all $i \in I$;

    \item $\sum_{i \in I} \phi_i \equiv 1$ on $O \subset \mathbb{R}^2$; and

    \item $\sup_{i \in I}
               \| \, \partial^\alpha \phi_i^\natural \, \|_{L^\infty}
            < \infty $ for all $\alpha \in \mathbb{N}_0^2$, where
            $\phi^\natural :
            \mathbb{R}^2 \to \mathbb{C},\xi\mapsto\phi_i(T_i \,\xi+b_i)$.
  \end{enumerate}
\label{eq:RPU}
\end{definition}

According to Theorem 2.4 in \cite{Bytchenkoff_2019}, there exists a regular partition of unity subordinate to the wave packet covering $\mathcal{Q}^{\alpha, \beta} = \{ Q_i \}_{i \in I}$ as the latter is almost structured.

\begin{definition}
The sequence $w = \{ w_i \}_{i \in I}$ of positive numbers is called $weight$. The weight is called $\mathcal{Q}$-$moderate$ where $\mathcal{Q} = \{ Q_i \}_{i \in I}$ stands for an almost structured covering of a Lebesgue-measurable subset $O$ of $\mathbb{R}^2$, if there exists such a positive number $C$ that $w_i \leqslant C \cdot w_{i'}$ for all such $i$ and $i' \in I$ that $Q_i \cap Q_{i'} \neq \emptyset$.
\label{eq:Gewicht}
\end{definition}

\begin{definition}
Let $0 \leqslant \beta \leqslant \alpha \leqslant 1$, $s \in \mathbb{R}$, $I_0$ be as defined by  \eqref{eq:Index-Menge}, and
  \begin{equation}
    w_i^{s}
    := \begin{cases}
         2^{j s} \, 
         & \text{if } i = (j, m, l) \in I_0 \setminus \{ (0, 0, 0) \} \,  \\
         1 \, 
         & \text{if } i = 0 =(0, 0, 0) \, 
       \end{cases}
  \hspace{0.25cm}.
  \label{eq:Wellenzuehen-Gewicht_}
  \end{equation}  
Then $w^{s} = (w_i^{s})_{i \in I_0}$ is called $wave$ $packet$ $weight$. 
\label{eq:Wellenzuehen-Gewicht}
\end{definition}  
  
The proof that $w^{s} = (w_i^{s})_{i \in I_0}$ is $\mathcal{Q}^{(\alpha,\beta)}$-moderate can be found in Lemma 6.1 of \cite{Bytchenkoff_2019}.

\begin{definition}
Let $\mathcal{Q} = (Q_i)_{i \in I_0}$, $\Phi = (\phi_i)_{i \in I_0}$,
  and $w = (w_i)_{i \in I_0}$ be an almost structured covering of a Lebesgue-measurable subset $O$ of $\mathbb{R}^2$,
  a regular partition of unity subordinate to $\mathcal{Q}$ and a $\mathcal{Q}$-moderate
  weight respectively and let $p,q \in (0,\infty]$. Then
  \[ \mathcal{D} (\mathcal{Q}, L^p,\ell_w^q)
    := \left\{ g \in Z' : \| g \|_{\mathcal{D} (\mathcal{Q}, L^p,\ell_w^q)} < \infty
       \right\} \,
  \] where $Z'$ stands for the topological dual of $Z := \mathcal{F} \left( C_c^\infty (\mathbb{R}^2) \right) \subset \mathcal{S}(\mathbb{R}^2)$ and
  $\| g \|_{\mathcal{D} (\mathcal{Q}, L^p,\ell_w^q)}$ for the quasi-norm
  \begin{equation}
    \| g \|_{\mathcal{D} (\mathcal{Q}, L^p,\ell_w^q)}
    := \left\|
          \left(
            w_i \cdot \| \mathcal{F}^{-1} (\varphi_i \cdot \widehat{g}) \|_{L^p}
          \right)_{i \in I_0}
       \right\|_{\ell^q} \in [0,\infty] \,
    \label{eq:DecompositionSpaceNorm}
  \end{equation} is called $decomposition$ $space$.
\label{eq:Decomposition_space}
\end{definition}

As explained in \cite{Bytchenkoff_2019}, thus defined decomposition space is a quasi-Banach space. In this definition, which differs slightly from that introduced originally \cite{Feichtinger_1985}, the functions $\{ \mathcal{F}^{-1} (\varphi_i \cdot \widehat{g})\}_{i \in I_0}$ play the role of components of the function $g$ localised in the frequency domain. Therefore, for \eqref{eq:DecompositionSpaceNorm} to be finite and so for the function $g$ to belong to $\mathcal{D} (\mathcal{Q}, L^p,\ell_w^q)$, these components must be $p$-integrable and their contributions $\{ \| \mathcal{F}^{-1} (\varphi_i \cdot \widehat{g}) \|_{L^p} \}_{i \in I_0}$ to the the norm $\| \bullet \|_{\mathcal{D} (\mathcal{Q}, L^p,\ell_w^q)}$ weighted with $w_i$ must be $q$-summable. With all these definitions we can now introduce the wave packet space.

\begin{definition}
Let $0 \leqslant \beta \leqslant \alpha \leqslant 1$, $s \in \mathbb{R}$, $p$ and $q \in (0, \infty ]$, and $\mathcal{Q}^{\alpha , \beta}$ and $w_i^{s}$ be the wave packet covering and weight respectively. The $wave$ $packet$ $space$ with parameters $\alpha,\beta, p, q$ and $s$ is defined as the decomposition space
\begin{equation*}
    \mathcal{W}_{s}^{p, q} (\alpha, \beta)
    := \mathcal{D} (\mathcal{Q}^{\alpha , \beta}, L^p,\ell_{w^s}^q) \, .
\end{equation*}
\label{eq:WPS}
\end{definition}

The parameter $s$ determines what weight $w_i^{s}$ will be assigned to the contribution of the $i$-th component $\mathcal{F}^{-1} (\varphi_i \cdot \widehat{g})$ of $g$ to the norm \eqref{eq:DecompositionSpaceNorm}. Therefore varying the value of $s$ amounts to making $\mathcal{W}_{s}^{p, q} (\alpha, \beta)$ more or, for that matter, less regular, since $i$ determines, among other things, the absolute value of the frequency at which the $i$-th component is localised. More details on the wave packet spaces, their properties and embeddings in each other or more classical function spaces such as Besov or Sobolev spaces can be found in \cite{Bytchenkoff_2019}.

\section{Banach frames and sets of atoms of wave packet spaces and their good localisation}

\noindent To be able to decompose or synthesise the elements of  $\mathcal{W}_{s}^{p, q} (\alpha, \beta)$ we shall need Banach frames and sets of atoms of it. Here are the definitions of these two notions \cite{Groechenig_1991}.

\begin{definition}
A set $\{\psi_i\}_{i \in I}$ in the dual space $X'$ of a quasi-Banach space $X$ is called $Banach$ $frame$ for $X$ if there exists well-defined bounded map, called $analysis$ $operator$, $A : X \rightarrow x , f \mapsto \{\langle \psi_i , f \rangle\}_{i \in I}$ where $x := \{\{\langle \psi_i , f \rangle\}_{i \in I} : f \in X\}$ is a solid Banach subspace of ${\mathbb{C}}^I$ and there exists such a bounded linear map $A_l^{-1} : x \mapsto X$ that $A_l^{-1} \circ A = I_X$ where $I_X$ stands for an identity operator on $X$.
\label{eq:Banach_frame}
\end{definition}

\begin{definition}
A set $\{\phi_i\}_{i \in I}$ in a quasi-Banach space $X$ is called $set$ $of$ $atoms$ in $X$ if there exists a well-defined bounded map, called $synthesis$ $operator$, $S : x \rightarrow X , \{c_i\}_{i \in I} \mapsto  \sum_{i \in I} c_i \phi_i$ where the coefficient space $x := \{c_i\}_{i \in I}$ associated with $\{\phi_i\}_{i \in I}$ is a solid subspace of ${\mathbb{C}}^I$ and there exists such a bounded linear map $S_r^{-1} : x \mapsto X$ that $S \circ S_r^{-1} = I_X$ where $I_X$ stands for an identity operator on $X$. The series expansion $g = \sum_{i \in I} c_i \phi_i$ of a given function $g \in X$ where $\{\phi_i\}_{i \in I}$ is a set of atoms is called $atomic$ $decomposition$ of $g$.
\label{eq:atomic_decomposition}
\end{definition}

In Definition 9.1 of \cite{Bytchenkoff_2019} we first introduced what we called $wave$ $packet$ $system$ and then, in Theorems 9.3 and 9.4 of the same report, established the conditions on the prototype functions of the elements of the system under which it will constitute either a Banach frame or a set of atoms for the wave packet space. By integrating these conditions with the definition, we shall now redefine the wave packet system so that it will automatically form both a Banach frame and a set of atoms of the corresponding wave packet space. 

\begin{definition}
First, let $0 \leqslant \beta \leqslant \alpha < 1$, $s_0 \leqslant 0$, $p_0$ and $q_0 \in (0,1]$. Second, let $t :=(t_1 , t_2) \in \mathbb{R}^2$, $\gamma (t)$ and $\psi (t) \in C^1 (\mathbb{R}^2)$, $\partial^a \gamma (t)$ and $\partial^a \psi (t) \in L^1 (\mathbb{R}^2) \cap L^{\infty} (\mathbb{R}^2)$ for any $a \in \mathbb{N}_0^2$ such that $\vert a \vert \leqslant 1$ and
\begin{equation}
\begin{split}
\sup_{t \in \mathbb{R}^2} (1 + |t|)^{\kappa_0} \cdot |\gamma (t)| < \infty \hspace{0.25cm} 
 \text{ and }
\hspace{0.25cm} \sup_{t \in \mathbb{R}^2} (1 + |t|)^{\kappa_0} \cdot |\psi (t)| < \infty
\end{split}
\label{eq:Bedingung_im Zeitgebiet}
\end{equation} where $\kappa_0 := 10 + 2 \cdot p_0^{-1}$. Third, let $\xi :=(\xi_1 , \xi_2) \in \mathbb{R}^2$, $\hat{\gamma} (\xi)$ and $\hat{\psi} (\xi) \in C^{\infty} (\mathbb{R}^2)$, $\hat{\gamma} (\xi) \neq 0$ for any $\xi \in \overline B_4 (0)$, $\hat{\psi} (\xi) \neq 0$ for any $\xi \in [-\epsilon, 1+\epsilon] \times [-1-\epsilon, 1+\epsilon]$,
\begin{equation}
\vert \partial^a \widehat{\partial^b \gamma (\xi)} \vert \lesssim \left( 1+\vert \xi \vert \right)^{-\kappa} \cdot \left( 1+\vert \xi_1 \vert \right)^{-\kappa_1} \cdot \left( 1+\vert \xi_2 \vert \right)^{-\kappa_2}
\label{eq:Bedingung_gamma}
\end{equation} and
\begin{equation}
\vert \partial^a \widehat{\partial^b \psi (\xi)} \vert \lesssim \left( 1+\vert \xi \vert \right)^{-\kappa} \cdot \left( 1+\vert \xi_1 \vert \right)^{-\kappa_1} \cdot \left( 1+\vert \xi_2 \vert \right)^{-\kappa_2}
\label{eq:Bedingung}
\end{equation} for any $a \in \mathbb{N}_0^2$ such that $\vert a \vert \leqslant 12$, any $b \in \mathbb{N}_0^2$ such that $\vert b \vert \leqslant 1$, $\kappa \geqslant 10$,  $\kappa_1 \geqslant 2$ and  $\kappa_2 \geqslant 10$. More precise, but rather cumbersome definitions of $a$, $b$, $\kappa$,  $\kappa_1$ and $\kappa_2$, which all depend on $\alpha$, $\beta$, $p_0$ and $q_0$, can be found in the statements of Theorems 9.3 and 9.4 of \cite{Bytchenkoff_2019}. Fourth, let $j^{\max} \in \mathbb{N}$,
\begin{equation}
I := \{ (0, 0, 0, k) : k \in \mathbb{Z}^2 \} \cup \{(j, m, l, k) \in \mathbb{N} \times \mathbb{N}_0 \times \mathbb{N}_0 \times \mathbb{Z}^2 : j \leqslant j^{\max} \text{, } m \leqslant m_j^{\max} \text{ and } l \leqslant l_j^{\max} \}
\label{eq:Index}
\end{equation} where $m_j^{\max}$ and $l_j^{\max}$ defined by \eqref{eq:m-Index-Grenzen} and \eqref{eq:l-Index-Grenzen} respectively and
\begin{equation} 
\psi_{i} (t) := \begin{cases}
 \gamma (t - \delta \cdot k) & \text{ if } (j, m, l)=(0, 0, 0) \\
 \vert  \det A_j \vert^{1/2} \cdot e^{-2\pi i \langle R_{jl} B_{jm}, t \rangle} \cdot \left( \psi \hspace{0.15cm} \circ A_{j}  \circ R_{jl}^{-1} \right) (t - \delta \cdot R_{jl} A_j^{-1} k)  & \text{ if } (j, m, l) \neq (0, 0, 0)
 \end{cases}
\label{eq:Wellenzug}
\end{equation} where $A_j$, $B_{jm}$ and $R_{jl}$ defined by \eqref{eq:A_j_und_B_jm} and \eqref{eq:R_jl}. Then the set $W (\alpha, \beta) := \{ \psi_i \}_{i \in I}$ is called $wave$ $packet$ $system$.
\label{eq:WPSy}
\end{definition}

Theorems 9.3 and 9.4 of \cite{Bytchenkoff_2019} indicate, as a corollary, that thus defined wave packets are both Banach frames and sets of atoms of $\mathcal{W}_{s}^{p\, q} (\alpha, \beta)$ with $s \in [-s_0, s_0]$, $p \in [p_0, \infty]$, $q \in [q_0, \infty]$ as long as $\delta \in (0, \delta_0]$ with $\delta_0 =\delta_0 (\alpha, \beta, s_0, p_0, q_0, \gamma, \psi) >0$ and that the coefficient space that is associated with $\mathcal{W}_{s}^{p\, q} (\alpha, \beta)$ through $\{ \psi_i \}_{i \in I}$ is given by the following definition.

\begin{definition}
Let $0 \leqslant \beta \leqslant \alpha < 1$, $s \in \mathbb{R}$, $p$ and $q \in (0,\infty]$, and $I_0$ be as defined in \eqref{eq:Index-Menge}. The set of sequences of complex numbers

\begin{equation*}
\mathcal{C}_{s}^{p \,q}(\alpha,\beta)
  := \left\{ c= \{ c_i^{(k)} \}_{i \in I_0, \, k \in \mathbb{Z}^2} \in \mathbb{C}^{I_0 \times \mathbb{Z}^2}
         \,:\,
         \|c\|_{\mathcal{C}_s^{p\,q}} < \infty
       \right\}
\end{equation*} where
\begin{equation*}
\|c\|_{\mathcal{C}_s^{p\,q}}
    := \left\|
         \left(
           w_i^{s + (\alpha + \beta) \cdot (\frac{1}{2} - \frac{1}{p})}
           \cdot \big\| (c_k^{(i)})_{k \in \mathbb{Z}^2} \big\|_{\ell^p}
         \right)_{i \in I_0}
       \right\|_{\ell^q}
    \in [0,\infty]
\end{equation*} is called $coefficient$ $space$ associated with the wave packet space $\mathcal{W}_{s}^{p\, q} (\alpha, \beta)$.
\label{eq:coefficient_space}
\end{definition}

The Fourier transform of $\psi_i (t)$ can be calculated to be
\begin{equation} 
\hat{\psi}_i (\xi) := \begin{cases}
 e^{-2\pi i \langle \delta \cdot k , \xi \rangle} \cdot \hat{\gamma} (\xi) & \text{ if } (j, m, l) = (0, 0, 0) \\
\begin{split} 
 \vert  \det A_j \vert^{1/2} & \cdot e^{-2\pi i \langle \delta \cdot R_{jl} A_j^{-1} k , \xi \rangle} \\ & \cdot \left( \hat{\psi} \circ A_{j}^{-1}  \circ R_{jl}^{-1} \right) (\xi - R_{jl} A_j B_{jm} )
\end{split} 
  & \text{ if } (j, m, l) \neq (0, 0, 0)
 \end{cases}
\hspace{0.25cm}.
\label{eq:FT_des_Wellezugs}
\end{equation}

The functions $\gamma (t)$ and $\psi (t)$ play the role of prototypes of the wave packets $\psi_i (t)$ with $i \in I$. Specifically all wave packets $\psi_i (t)$ with $i$ such that $(j, m, l) = (0, 0, 0)$ are generated from $\gamma (t)$ only by shifting it, while the wave packets $\psi_i (t)$ with $i$ such that $(j, m, l) \neq (0, 0, 0)$ are generated from $\psi (t)$ by scaling, modulating, rotating and shifting it. The index $j \in \mathbb{N}$ identifies the measurements $2^{-\alpha j} \times 2^{-\beta j}$ of the wave packet $\psi_i (t)$ -- compared with those of its prototype -- and, indeed,  the measurements $2^{\alpha j} \times 2^{\beta j}$ of its Fourier transform $\hat{\psi}_i (\xi)$ -- compared with those of the Fourier transform of its prototype. The indexes $j$, $m$ and $l$ determine the position $\xi_i$ of the maximum of the Fourier transform $\hat{\psi}_i (\xi)$ of the wave packet $\psi_i (t)$ in the frequency plane, namely

\begin{equation} 
\xi_i := \xi_{jml} := \begin{cases}
 \left( \begin{array}{ cc }
0 \\
0  \\
\end{array} \right) & \text{ if } (j, m, l) = (0, 0, 0) \\
 R_{jl} B_{jm}
=  R_{jl}
\left( \begin{array}{ cc }
2^{j-1} + m \cdot 2^{\alpha j} \\
0  \\
\end{array} \right) & \text{ if } (j, m, l) \neq (0, 0, 0)
 \end{cases}
\hspace{0.25cm}.
\label{eq:Ort_im_Frequenzbereich}
\end{equation} The indexes $j$ and $m$ in particular determine the distance $r_i$ of the maximum of the of Fourier transform $\hat{\psi}_i (\xi)$ of the wave packet $\psi_i (t)$ from the origin of the frequency plane, namely
\begin{equation} 
r_i := \vert \xi_i \vert := \vert \xi_{jml} \vert := \begin{cases}
 0 & \text{ if } (j, m, l) = (0, 0, 0) \\
 2^{j-1} + m \cdot 2^{\alpha j} & \text{ if } (j, m, l) \neq (0, 0, 0)
 \end{cases}
\hspace{0.25cm},
\label{eq:Abstand_in_Frequenzbereich}
\end{equation}
while the indexes $j$ and $l$ determine the angle $\theta_i$ at which the Fourier transform  $\hat{\psi}_i (\xi)$ of the wave packet $\psi_i (t)$ is inclined to the $\xi_1$-axis in the frequency plane, namely
\begin{equation} 
\theta_i := \begin{cases}
 0 & \text{ if } (j, m, l) = (0, 0, 0) \\
 \theta_{jl} & \text{ if } (j, m, l) \neq (0, 0, 0)
 \end{cases}
\label{eq:Winckel}
\end{equation} with $\theta_{jl}$ defined by \eqref{eq:theta_jl}. In other words $r_i$ and $\theta_i$ are polar coordinates of the maximum of  $\hat{\psi}_i (\xi)$ in the frequency plane. The indexes $k$, $j$ and $l$ determine the position of the maximum of the wave packet $\psi_i (t)$ in the time plane, namely

\begin{equation} 
t_i:= t_k^{(jl)} := \begin{cases}
 \delta \cdot k =
\delta \cdot
\left( \begin{array}{ cc }
k_1 \\
k_2 \\
\end{array} \right) & \text{ if } (j, m, l) = (0, 0, 0) \\
 \delta \cdot R_{jl} A_{j}^{-1} k
= \delta \cdot R_{jl}
\left( \begin{array}{ cc }
k_1 \cdot 2^{-\alpha j} \\
k_2  \cdot 2^{-\beta j}  \\
\end{array} \right) & \text{ if } (j, m, l) \neq (0, 0, 0)
 \end{cases}
\hspace{0.25cm},
\label{eq:Ort_im_Zeitbereich}
\end{equation}

Whether the set $\{\psi_i\}_{i \in I}$ of the wave packets can be orthogonalised, say, by imposing additional conditions on the prototypes $\gamma$ and $\phi$ is the subject of a separate ongoing project. At the moment we shall show that $\langle \psi_i \vert \psi_{i'} \rangle$ tends to zero as $i$ and $i'$ become different. To measure the difference between $i$ and $i'$ we shall need the following definition.

\begin{definition} Let $i=(j, m, l, k)$ and $i'=(j', m', l', k') \in I$ where $I$ defined by \eqref{eq:Index}, $r (i, i') := \vert r_i -r_{i'} \vert \in \mathbb{R}_+$ where $r_i$ and $r_{i'}$ defined by \eqref{eq:Abstand_in_Frequenzbereich}, $\theta (i, i') := \vert \theta_i -\theta_{i'} \vert \in [0, \pi]$ where $\theta_i$ and $\theta_{i'}$ defined by \eqref{eq:Winckel}, $ t_1 (i,i') := \vert (t_i-t_{i'})_1  \vert  \in \mathbb{R}_+$ and $ t_2 (i,i') := \vert (t_i-t_{i'})_2 \vert \in \mathbb{R}_+$ where $t_i$ and $t_{i'}$ defined by \eqref{eq:Ort_im_Zeitbereich}. Then
\begin{equation}
\rho ( i, i') := r ( i, i') + \theta ( i, i') + t_1  ( i, i') + t_2 (i, i')
\label{eq:rho}
\end{equation} will be called $distance$ between indexes $i$ and $i'$ of the wave packets $\psi_i$ and  $\psi_{i'}$.
\end{definition}

Thus the distance $\rho( i , i' )$ between the indexes $i$ and $i'$ and, indeed, between the wave packets $\psi_i$ and $\psi_{i'}$ indexed by them is defined as the sum of the lengths $\vert (t_i-t_{i'})_1  \vert$ and $\vert (t_i-t_{i'})_2 \vert$ of the orthogonal projections of the difference $t_i - t_{i'}$ between the radii vectors $t_i$ and $t_{i'}$ from the origin of the time plane to the maxima of the wave packets $\psi_i$ and $\psi_{i'}$ on the horizontal and vertical axes of the Cartesian coordinate system of the time plane and the absolute values of the differences $r ( i, i')$ and $\theta ( i, i')$ between the polar coordinates of the maxima of the Fourier transforms $\hat{\psi}_i$ and $\hat{\psi}_{i'}$ of the wave packets $\psi_i$ and $\psi_{i'}$ on the frequency plane. The following Lemma justifies the name $distance$ of $\rho( i, i')$.
\newline

\noindent \textbf{Lemma 1.} Let $j^{\max}$ be a finite natural number and $i, i' \in I$ where $I$ defined by \eqref{eq:Index}, then
\begin{equation}
\rho( i , i' ) \geqslant 0
\hspace{0.25cm},
\label{eq:Nicht-Negativität}
\end{equation}
\begin{equation}
i=i' \Leftrightarrow \rho( i , i' ) = 0
\hspace{0.25cm},
\label{eq:Indiscernabiliy}
\end{equation}
\begin{equation}
\rho( i , i' ) = \rho( i' , i )
\label{eq:Symmetrie}
\end{equation} and
\begin{equation}
\rho( i , i'' ) \leqslant \rho( i , i' ) + \rho( i' , i'' )
\hspace{0.25cm}.
\label{eq:Dreieck}
\end{equation}
\noindent $Proof$. First, the non-negativity of $\rho (i', i)$ expressed by \eqref{eq:Nicht-Negativität} follows from its being the sum of four non-negative numbers.

Second, if $i=(j, m, l, k) = i'=(j', m', l', k')$, then $r (i, i')$, $\theta (i, i')$, $t_1 (i, i')$ and $t_2 (i, i')$  equal zero according to \eqref{eq:Abstand_in_Frequenzbereich},  \eqref{eq:Winckel} and \eqref{eq:Ort_im_Zeitbereich} respectively. This, in its turn, implies that $\rho  (i , i') = 0$. If, conversely, $\rho  (i , i') = 0$, then $\theta (i, i') = 0$. These two equalities imply that $j = j'$. Indeed, if $j'$ equalled, say, $j+1$, then, on the one hand, $r_{i'}$ would equal or be larger than $2^{j'-1} = 2^{j}$, according to \eqref{eq:Abstand_in_Frequenzbereich}, and, on the other hand,
\begin{equation*}
r_{i} \leqslant 2^{j-1} + m_j^{\max} \cdot 2^{\alpha j}
  =  2^{j-1} + (\lceil 2^{(1-\alpha) j -1} \rceil - 1) \cdot 2^{\alpha j}
  < 2^{j-1} + 2^{(1-\alpha) j -1} \cdot 2^{\alpha j} 
  = 2^j
\hspace{0.25cm},
\end{equation*} according to \eqref{eq:Abstand_in_Frequenzbereich} and \eqref{eq:m-Index-Grenzen}. This, in its turn, would imply that $\vert r_i - r_{i'} \vert > 0$, which would contradict $\rho  (i , i') = 0$. If, otherwise, $j'$ equalled or were larger than $j+2$, then, on the one hand, $r_{i'}$ would equal or be larger than $2^{j'-1} = 2^{j+1}$, according to \eqref{eq:Abstand_in_Frequenzbereich}, and, on the other hand,
\begin{equation*}
r_{i} \leqslant 2^{j-1} + m_j^{\max} \cdot 2^{\alpha j}
  =  2^{j-1} + \lceil 2^{(1-\alpha) j -1} \rceil \cdot 2^{\alpha j}
  < 2^{j-1} + \left( 2^{(1-\alpha) j -1} +1 \right) \cdot 2^{\alpha j} 
  = 2^j + 2^{\alpha j} \leqslant 2^{j+1}
\hspace{0.25cm},
\end{equation*} according to \eqref{eq:Abstand_in_Frequenzbereich} and \eqref{eq:m-Index-Grenzen}. This would again imply that $\vert r_i - r_{i'} \vert > 0$, which would contradict $\rho  (i , i') = 0$. In a similar manner it can be demonstrated that $j' < j$ and $\rho (i, i') = 0$ are not compatible either. Furthermore, according to \eqref{eq:Winckel}, \eqref{eq:theta_jl} and \eqref{eq:l-Index-Grenzen}, $0 \leqslant \theta_i < 2 \pi$ and so $0 \leqslant \theta (i, i') < 2 \pi$ too. Therefore either $l = 0 = l'$ if $i$ is such that $j = j' = 0$ or
\begin{equation*}
\theta (i, i') = 2 \pi \cdot \frac{{2}^{ \left( \beta -1 \right) j}}{N} \cdot \vert l - l' \vert = 0
\hspace{0.25cm}
\end{equation*} if it is not, which, given ${2}^{ \left( \beta -1 \right) j} > 0$ for any finite $j \leqslant j^{\max}$, again indicates that $l=l'$. Moreover $\rho (i, i') = 0$ impies that $r (i, i') = 0$ or, according to \eqref{eq:Abstand_in_Frequenzbereich}, that $2^{j-1} + m \cdot 2^{\alpha j} = 2^{j'-1} + m' \cdot 2^{\alpha j'}$. This, together with $j=j'$, implies that $m=m'$. Finally, according to \eqref{eq:Ort_im_Zeitbereich},
\begin{equation} 
t_i - t_{i'}
   = \delta \cdot R_{jl}
     \left( \begin{array}{ cc }
               (k_1 - k'_{1}) \cdot 2^{-\alpha j} \\
               (k_2 - k'_{2}) \cdot 2^{-\beta j}  \\
             \end{array} \right)
\label{eq:t_i-t_i'} 
\end{equation} as $j=j'$ and $l=l'$ and so $R_{jl} = R_{j'l'}$. Therefore $\rho (i, i') = 0$ implies that
\begin{equation*} 
\begin{cases}
  (t_i - i_{i'})_1
   = \delta \cdot \cos \theta_{jl} \cdot (k_1 - k'_{1}) \cdot 2^{-\alpha j}
    - \delta \cdot \sin \theta_{jl} \cdot (k_2 - k'_{2}) \cdot 2^{-\beta j}
   = 0  \\
  (t_i - i_{i'})_2
   = \delta \cdot \sin \theta_{jl} \cdot (k_1 - k'_{1}) \cdot 2^{-\alpha j}
    + \delta \cdot \cos \theta_{jl} \cdot (k_2 - k'_{2}) \cdot 2^{-\beta j}
   = 0
 \end{cases}
\hspace{0.25cm}.
\label{eq:lineare_Gleichung}
\end{equation*} The determinant of this homogeneous system of two linear equations with unknowns $k_1 - k'_{2}$ and $k_1 - k'_{2}$ equals $\delta^2 \cdot 2^{-(\alpha + \beta) j}$, which is a positive number for any finite $j \leqslant j^{\max}$. Therefore the system has the trivial solution $k_1 - k'_{2} =0$ and $k_1 - k'_{2} =0$ only. This completes the proof of the fact that $\rho (i, i') =0$ implies $i=i'$.

Third, the symmetry of $\rho (i, i')$ expressed by \eqref{eq:Symmetrie} results from that of $r (i, i')$, $\theta (i, i')$, $t_1 (i, i')$ and $t_2 (i, i')$. Fourth, the subadditivity of $\rho( i, i')$ expressed by \eqref{eq:Dreieck} results from that of $r (i, i')$, $\theta (i, i')$, which are the absolute values of the differences between real numbers, and that of $t_1 (i, i')$ and $t_2 (i, i')$, which are the absolute values of the orthogonal projections of the difference between two vectors. $\Box$
\newline

The next lemma will allow us to obtain one of the major results of this work and is a corollary to Lemma 1.
\newline

\noindent \textbf{Lemma 2.} Let $j^{\max}$ be a finite natural number, then the metric index space $I$ defined by \eqref{eq:Index} with the distance defined by \eqref{eq:rho} is separated, i.e.
\begin{equation}
\inf_{i, \, i' \in I, \, i \neq i'} \rho ( i, i') = C > 0
\hspace{0.25cm}.
\label{eq:getrennte_Index-Menge}
\end{equation}

\noindent $Proof$. We note that
\begin{equation*}
\inf_{i, \, i' \in I, \, i \neq i'} \rho ( i, i')
\geqslant \inf_{i, \, i' \in I, \, i \neq i'} \left( r ( i, i') + \theta ( i, i') \right)
+ \inf_{i, \, i' \in I, \, i \neq i'}  \left( t_1 ( i, i') + t_2 ( i, i') \right)
\label{eq:Infimum}
\end{equation*} and consider two possible situations that can occur. First let $i$ and $i'$ be such that $j \neq j'$ and/or $m \neq m'$ and/or $l \neq l'$. From Lemma 1 we know that $\rho ( i, i') > 0$ for any $i$ and $i'$ that differ one way or another from each other. In particular $r ( i, i') + \theta ( i, i') = \rho ( i, i') > 0$ as $k_1 = k'_1$ and $k_2 = k'_2$ and $j \neq j'$ and/or $m \neq m'$ and/or $l \neq l'$. From this and the fact that $r ( i, i') + \theta ( i, i')$ does not depend on $k_1$, $k_2$, $k'_1$ and $k'_2$ we infer that $r ( i, i') + \theta ( i, i') > 0$ as $j \neq j'$ and/or $m \neq m'$ and/or $l \neq l'$, no matter what $k_1$, $k_2$, $k'_1$ and $k'_2$ are. As $j^{\max}$ is finite, so are $m^{\max}$ and $l^{\max}$. Therefore, for given $k_1$, $k_2$, $k'_1$ and $k'_2$, there is only a finite number of different pairs of indexes $(i, i')$ -- in which $j \neq j'$ and/or $m \neq m'$ and/or $l \neq l'$ -- and amongst them there is a pair that minimises $r ( i, i') + \theta ( i, i')$. Let us call this positive minimum $C_1$, then
\begin{equation}
\inf_{i, \, i' \in I, \, i \neq i'} \rho ( i, i')
\geqslant \inf_{i, \, i' \in I, \, i \neq i'} \left( r ( i, i') + \theta ( i, i') \right)
> C_1 > 0 
\label{eq:Infimum_1}
\end{equation} as $i$ and $i'$ are such that $j \neq j'$ and/or $m \neq m'$ and/or $l \neq l'$.
 
Let now $i$ and $i'$ be such that $j = j'$, $m = m'$ and $l = l'$. Then $R_{jl} = R_{j'l'}$ and $t_i-t_{i'}$ is given by \eqref{eq:t_i-t_i'}, from which we infer that 
\begin{equation*} 
\begin{split}
t_1 (i, i') + t_2 (i, i') & = \vert (t_i - i_{i'})_1 \vert + \vert (t_i - i_{i'})_2 \vert
 \geqslant \vert \cos \theta_{jl} \vert \cdot \vert (t_i - i_{i'})_1 \vert
     + \vert \sin \theta_{jl} \vert \cdot \vert (t_i - i_{i'})_2 \vert \\ 
 & = \vert \cos \theta_{jl} \cdot (t_i - i_{i'})_1 \vert
     + \vert \sin \theta_{jl} \cdot (t_i - i_{i'})_2 \vert \\
 & \geqslant \vert \cos \theta_{jl} \cdot (t_i - i_{i'})_1
     + \sin \theta_{jl} \cdot (t_i - i_{i'})_2 \vert \\   
 &  = \big\vert \delta \cdot \cos^2 \theta_{jl} \cdot (k_1 - k'_{1}) \cdot 2^{-\alpha j}
     - \delta \cdot \cos \theta_{jl} \cdot \sin \theta_{jl} \cdot (k_2 - k'_{2}) \cdot 2^{-\beta j} \\
 & + \delta \cdot \sin^2 \theta_{jl} \cdot (k_1 - k'_{1}) \cdot 2^{-\alpha j}
     + \delta \cdot \sin \theta_{jl} \cdot \cos \theta_{jl} \cdot (k_2 - k'_{2}) \cdot 2^{-\beta j} \big\vert \\
 & = \delta \cdot \vert k_1 - k'_{1} \vert \cdot 2^{-\alpha j}    
\end{split}
\end{equation*} and
\begin{equation*} 
\begin{split}
t_1 (i, i') + t_2 (i, i') & \geqslant \vert \sin \theta_{jl} \vert \cdot \vert (t_i - i_{i'})_1 \vert
     + \vert \cos \theta_{jl} \vert \cdot \vert (t_i - i_{i'})_2 \vert \\ 
 & = \vert - \sin \theta_{jl} \cdot (t_i - i_{i'})_1 \vert
     + \vert \cos \theta_{jl} \cdot (t_i - i_{i'})_2 \vert \\
 & \geqslant \vert -\sin \theta_{jl} \cdot (t_i - i_{i'})_1
     + \cos \theta_{jl} \cdot (t_i - i_{i'})_2 \vert \\   
 &  = \big\vert - \delta \cdot \sin \theta_{jl} \cdot \cos \theta_{jl} \cdot (k_1 - k'_{1}) \cdot 2^{-\alpha j}
     + \delta \cdot \sin^2 \theta_{jl} \cdot (k_2 - k'_{2}) \cdot 2^{-\beta j} \\
 &   + \delta \cdot \cos \theta_{jl} \cdot \sin \theta_{jl} \cdot (k_1 - k'_{1}) \cdot 2^{-\alpha j}
     + \delta \cdot \cos^2 \theta_{jl} \cdot (k_2 - k'_{2}) \cdot 2^{-\beta j} \big\vert \\
 & = \delta \cdot \vert k_2 - k'_2 \vert \cdot 2^{-\beta j}    
\end{split}
\hspace{0.25cm}.
\end{equation*} Thus
\begin{equation*}
\begin{split}
t_1 (i, i') + t_2 (i, i') &
  \geqslant 1/2 \cdot \left( \delta \cdot \vert k_1 - k'_1 \vert \cdot 2^{-\alpha j}
  + \delta \cdot \vert k_2 - k'_2 \vert \cdot 2^{-\beta j} \right)
  \geqslant 1/2 \cdot \delta \cdot 2^{-\alpha j}
  \geqslant \delta \cdot 2^{-(\alpha j^{\max} +1)} >0
\end{split}
\end{equation*} if $k_1 \neq k'_1$ and/or $k_2 \neq k'_2$. Therefore
\begin{equation}
\inf_{i, \, i' \in I, \, i \neq i'} \rho ( i, i')
\geqslant \inf_{i, \, i' \in I, \, i \neq i'}  \left( t_1 ( i, i') + t_2 ( i, i') \right)
> \delta \cdot 2^{-(\alpha j^{\max} +1)} > 0
\label{eq:Infimum_2}
\end{equation} as $k_1 \neq k'_1$ and/or $k_2 \neq k'_2$ and $j = j'$, $m = m'$ and $l = l'$. Combining \eqref{eq:Infimum_1} and \eqref{eq:Infimum_2} results in \eqref{eq:getrennte_Index-Menge} where $C=\min (C_1, \,  \delta \cdot 2^{-(\alpha j^{\max} +1)})$. $\Box$
\newline

The next Lemma is an immediate corollary of the previous one and Lemma 3 in \cite{Groechenig_2003}.
\newline

\noindent \textbf{Lemma 3.} Let $j^{\max}$ and $I$ be a finite natural number and the set $I$ defined by \eqref{eq:Index} respectively. Then the distance defined by \eqref{eq:rho} satisfies the following inequality
\begin{equation}
\sup_{i \in I} \, \sum\limits_{i'} \left(1+ \rho ( i, i') \right)^{-n} < \infty
\label{eq:Supreimum}
\end{equation} for any $n>5$.
\newline

We now state and prove the main theorem of this section.
\newline

\noindent \textbf{Theorem 1.} Let $j^{\max}$ be a finite natural number, $\kappa_0$ and $\kappa$ in Definition 10 not smaller than 24, $\psi_i(t)$ and $\psi_{i'}(t)$ any two elements of the wave packet system $W (\alpha, \beta)$ of the wave packet space $\mathcal{W}_{s}^{p, q} (\alpha, \beta)$ and $\rho( i, i')$ the distance between their indexes $i$ and $i' \in I$ where $I$ defined by \eqref{eq:Index}, then
\begin{equation}
\vert \langle \psi_{i} \vert \psi_{i'} \rangle \vert  \lesssim  \left(1+ \rho(i, i') \right)^{-6}
\hspace{0.25cm}.
\label{eq:Produkt}
\end{equation}

\noindent $Proof$. From \eqref{eq:Bedingung_im Zeitgebiet} we infer that
\begin{equation}
\begin{split}
|\gamma (t)| \lesssim \frac{1}{(1 + |t|)^{\kappa_0}}
\hspace{0.25cm} \text{ and } \hspace{0.25cm}
|\psi (t)| \lesssim \frac{1}{(1 + |t|)^{\kappa_0}} 
\end{split}
\label{eq:Abschaetzung_von_Gamma}
\end{equation} and note that the right term in these inequalities is isotropic in the sense that it does not depend on the orientation of $t$ on the time plane. Therefore the functions obtained from $\gamma$ or $\psi$ by rotating them around the origin of the time plane satisfy the same inequalities. We also note that $1 \leqslant \det A_j = 2^{(\alpha + \beta) j} \leqslant 2^{(\alpha + \beta) j^{\max}} < \infty$ as long as $j^{\max}$ is finite and that $\vert\vert A_j \vert\vert \geqslant 1$. Furthermore from Section K.1 in \cite{Grafakos_2008} we know that
\begin{equation}
\int\limits_{\mathbb{R}^n}
\frac{1}{\left( 1+\vert x - a \vert \right)^m} \cdot \frac{1}{\left( 1+\vert x - b \vert  \right)^m} \hspace{0.1cm} dx \lesssim \frac{1}{\left( 1+ \vert a - b \vert \right)^m}
\label{eq:Grafakos}
\end{equation} for $m>n$. Therefore
\begin{equation}
\begin{split}
\vert \left\langle \psi_i \vert \psi_{i'} \right\rangle \vert &
   \leqslant \int\limits_{\mathbb{R}^2}
\vert \psi_i \vert \cdot \vert \psi_{i'} \vert 
 \hspace{0.1cm} dt
  \lesssim \int\limits_{\mathbb{R}^2}
\frac{1}{\left( 1+\vert t - t_i \vert \right)^{\kappa_0}} \cdot \frac{1}{\left( 1+\vert t - t_{i'} \vert  \right)^{\kappa_0}} \hspace{0.1cm} dt
  \lesssim \frac{1}{\left( 1+ \vert t_i - t_{i'} \vert \right)^{\kappa_0}} \\
 &  \leqslant \frac{1}{\left( 1+ \vert ( t_i - t_{i'} )_1 \vert \right)^{\frac{\kappa_0}{2}}}
       \cdot \frac{1}{\left( 1+ \vert ( t_i - t_{i'} )_2 \vert \right)^{\frac{\kappa_0}{2}}}
    \leqslant \frac{1}{\left( 1+ \vert ( t_i - t_{i'} )_1 \vert + \vert ( t_i - t_{i'} )_2 \vert \right)^{\frac{\kappa_0}{2}}} \\  
 &   = \frac{1}{\left( 1+ t_1  ( i, i') + t_2 (i, i') \right)^{\frac{\kappa_0}{2}}}
\end{split}
\hspace{0.25cm}.
\label{eq:Grafakos_im_Zeitgebiet}
\end{equation}

From \eqref{eq:Bedingung_gamma} and \eqref{eq:Bedingung} we infer that
\begin{equation}
\vert \hat{\gamma} (\xi) \vert \lesssim \frac{1}{\left( 1+\vert \xi \vert \right)^{\kappa}}
\hspace{0.25cm} \text{ and } \hspace{0.25cm}
\vert \hat{\psi} (\xi) \vert \lesssim \frac{1}{\left( 1+\vert \xi \vert \right)^{\kappa}}
\hspace{0.25cm}.
\label{eq:Bedingung_isotropisch}
\end{equation} From this and the arguments similar to those used to obtain \eqref{eq:Grafakos_im_Zeitgebiet} we infer that
\begin{equation}
\begin{split}
\left \vert \left\langle \hat{\psi}_i \big\vert \hat{\psi}_{i'} \right\rangle \right \vert &
   \leqslant  \int\limits_{\mathbb{R}^2}
\big\vert \hat{\psi}_i \big\vert \cdot \big\vert \hat{\psi}_{i'} \big\vert \hspace{0.1cm} d\xi
   \lesssim \int\limits_{\mathbb{R}^2}
   \frac{1}{\left( 1+\vert \xi - \xi_i \vert \right)^\kappa}
   \cdot \frac{1}{\left( 1+\vert \xi - \xi_{i'} \vert \right)^\kappa} \hspace{0.1cm} d\xi \\ &
   \lesssim \frac{1}{\left( 1+\vert \xi_i - \xi_{i'} \vert \right)^\kappa}
\end{split}
\hspace{0.25cm}.
\label{eq:Erste_Abschaetzung_0}
\end{equation} We now note that
\begin{equation}
\begin{split}
\frac{1}{1+ \vert \xi_i - \xi_{i'} \vert} & = 
\frac{1}{1+ \vert \xi_i - \xi_{i'} \vert \cdot \vert e^{-i\theta_{i'}} \vert} 
   = \frac{1}{1+ \vert \xi_i \cdot e^{-i\theta_{i'}} - \xi_{i'} \cdot e^{-i\theta_{i'}} \vert} \\ &
   = \frac{1}{1+ \left\vert \xi_i \cdot e^{-i\theta_{i'}} - r_{i'} \right\vert}
   \leqslant \frac{1}{1+ \left\vert \vert \xi_i \cdot e^{-i\theta_{i'}} \vert - r_{i'} \right\vert}
   = \frac{1}{1+ \vert r_i - r_{i'} \vert}
   = \frac{1}{1 + r (i, i')}
\end{split}
\label{eq:Abschaetzung_1}
\end{equation} for any $i$ and $i' \in I$. Furthermore, if $i$ is such that $(j, m, l) \neq (0, 0, 0)$ and so $r_{i} = 2^{j-1} + m \cdot 2^{\alpha j'} \geqslant 1$, then, on the one hand,
\begin{equation}
\begin{split}
\frac{1}{1+ \vert \xi_i - \xi_{i'} \vert} & = 
\frac{1}{1+ \vert \xi_i - \xi_{i'} \vert \cdot \vert e^{-i\theta_{i'}} \vert} 
   = \frac{1}{1+ \vert \xi_i \cdot e^{-i\theta_{i'}} - \xi_{i'} \cdot e^{-i\theta_{i'}} \vert} \\ &
   \leqslant \frac{1}{1+ \left\vert \left( \xi_i \cdot e^{-i\theta_{i'}} - \xi_{i'} \cdot e^{-i\theta_{i'}} \right)_2 \right\vert}
   = \frac{1}{1+ r_i \cdot \vert \sin (\theta_i - \theta_{i'}) \vert}
   \leqslant \frac{1}{1+ \vert \sin (\theta_i - \theta_{i'}) \vert} \\ &
   =  \frac{1}{1+ \sin ( \vert \theta_i - \theta_{i'} \vert )}
   \leqslant \frac{1}{1+ \frac{2}{\pi} \cdot \vert \theta_i - \theta_{i'} \vert}
   \leqslant \frac{\pi}{2} \cdot \frac{1}{1 + \vert \theta_i - \theta_{i'} \vert}
\end{split}
\label{eq:Abschaetzung_2.1}
\end{equation} as $\vert \theta_i - \theta_{i'} \vert \in [0, \frac{\pi}{2}]$; and, on the other hand,
\begin{equation}
\begin{split}
\frac{1}{1+ \vert \xi_i - \xi_{i'} \vert} = & 
\frac{1}{1+ \vert \xi_i - \xi_{i'} \vert \cdot \vert e^{-i (\theta_{i'} + \pi/2)} \vert} 
   = \frac{1}{1+ \left\vert \xi_i \cdot e^{-i (\theta_{i'} + \pi/2)} - \xi_{i'} \cdot e^{-i (\theta_{i'} + \pi/2)} \right\vert} \\ &
   \leqslant \frac{1}{1+ \left\vert \left( \xi_i \cdot e^{-i (\theta_{i'} + \pi/2)} - \xi_{i'} \cdot e^{-i (\theta_{i'} + \pi/2)} \right)_1 \right\vert} \\ &
   = \frac{1}{1+ r_i \cdot \vert \cos (\theta_i - \theta_{i'} -\pi/2) \vert}
   \leqslant \frac{1}{1+ \vert \cos (\theta_i - \theta_{i'} -\pi/2) \vert} \\ &
   = \frac{1}{1+ \vert \sin (\theta_i - \theta_{i'} - \pi) \vert} 
   \leqslant \frac{1}{1+ \frac{2}{\pi} \cdot \vert \theta_i - \theta_{i'} - \pi \vert} 
   \lesssim \frac{1}{1+ \vert \theta_i - \theta_{i'} \vert}     
\end{split}
\label{eq:Abschaetzung_2.2}
\end{equation} as $\vert \theta_i - \theta_{i'} \vert \in [ \frac{\pi}{2} , \pi]$ and so $\vert \theta_i - \theta_{i'} -\pi \vert \in [0, \frac{\pi}{2}]$. Indeed, if $\theta_i - \theta_{i'} \in [ \frac{\pi}{2} , \pi]$, then $\theta_i - \theta_{i'} -\pi \in [0, - \frac{\pi}{2}]$ and
\begin{equation}
\begin{split}
\frac{1}{1+ \frac{2}{\pi} \cdot \vert \theta_i - \theta_{i'} - \pi \vert} &
   \leqslant \frac{1}{ \left\vert 1- \frac{2}{\pi} \cdot \vert \theta_i - \theta_{i'} - \pi \vert \right\vert}
   = \frac{1}{ \left\vert 1 - \frac{2}{\pi} \cdot \left( \pi - \vert \theta_i - \theta_{i'} \vert \right) \right\vert } \\ &
   = \frac{1}{ \left\vert 1 + \frac{2}{\pi} \cdot \vert \theta_i - \theta_{i'} \vert \right\vert } 
   \leqslant \frac{\pi}{2} \cdot \frac{1}{1+ \vert \theta_i - \theta_{i'} \vert} 
\end{split}
\hspace{0.25cm};
\label{eq:Abschaetzung_2.2.1}
\end{equation} and, if $\theta_i - \theta_{i'} \in [ - \frac{\pi}{2} , - \pi]$, then $\theta_i - \theta_{i'} -\pi \in [0, \frac{\pi}{2}]$ and
\begin{equation}
\begin{split}
\frac{1}{1+ \frac{2}{\pi} \cdot \vert \theta_i - \theta_{i'} - \pi \vert} &
   \leqslant \frac{3}{3+ \frac{2}{\pi} \cdot \vert \theta_i - \theta_{i'} - \pi \vert}
   = \frac{3}{3+ \frac{2}{\pi} \cdot \left( \vert \theta_i - \theta_{i'} \vert - \pi \right)} \\ &
   = \frac{3}{1+ \frac{2}{\pi} \cdot \vert \theta_i - \theta_{i'} \vert} 
   \leqslant \frac{3 \pi}{2} \cdot \frac{1}{1+ \vert \theta_i - \theta_{i'} \vert}  
\end{split}
\hspace{0.25cm}.
\label{eq:Abschaetzung_2.2.2}
\end{equation}

From \eqref{eq:Abschaetzung_2.1}, \eqref{eq:Abschaetzung_2.2} and \eqref{eq:Abschaetzung_1} we infer that
\begin{equation}
\begin{split}
\frac{1}{1+ \vert \xi_i - \xi_{i'} \vert} 
  \lesssim \frac{1}{1+ \vert \theta_i - \theta_{i'} \vert}
  = \frac{1}{1+ \theta (i, i')}
\end{split}
\label{eq:Abschaetzung_2}
\end{equation} and
\begin{equation}
\begin{split}
\frac{1}{\left( 1 + \vert \xi_i - \xi_{i'} \vert \right)^2} 
  \lesssim \frac{1}{1+ r (i, i')} \cdot \frac{1}{1+ \theta (i, i')}
  \leqslant \frac{1}{1 + r (i, i') + \theta (i, i')}
\end{split}
\label{eq:Abschaetzung_2_endgueltig}
\end{equation} for any such $i \in I$ that $(j, m, l) \neq (0, 0, 0)$ and/or  $(j', m', l') \neq (0, 0, 0)$.

We now consider two possible situations. First, let $i$ and $i'$ be such that $(j, m, l) = (0, 0, 0) = (j', m', l')$. Then $r_i = r_{i'} = \theta_i = \theta_{i'} = 0$ and so $r (i, i') = \theta (i, i') = 0$ and $\rho (i, i') = t_1  ( i, i') + t_2 (i, i')$. Moreover $\kappa_0 \geqslant 12$ according to Definition 10. Therefore \eqref{eq:Grafakos_im_Zeitgebiet} implies \eqref{eq:Produkt} for any such $i$ and $i'$ that $(j, m, l) = (0, 0, 0) = (j', m', l')$. Second, if $i$ and $i'$ are such that $(j, m, l) \neq (0, 0, 0)$ and/or  $(j', m', l') \neq (0, 0, 0)$, then from \eqref{eq:Grafakos_im_Zeitgebiet}, \eqref{eq:Erste_Abschaetzung_0} and \eqref{eq:Abschaetzung_2_endgueltig} and the fact that $\vert \left\langle \psi_i \vert \psi_{i'} \right\rangle \vert = \left \vert \left\langle \hat{\psi}_i \big\vert \hat{\psi}_{i'} \right\rangle \right \vert$ we infer that
\begin{equation}
\begin{split}
\vert \left\langle \psi_i \vert \psi_{i'} \right\rangle \vert &
    = \left\vert \vert \left\langle \psi_i \vert \psi_{i'} \right\rangle \vert^{\frac{1}{2}}
    \cdot \left \vert \left\langle \hat{\psi}_i \big\vert \hat{\psi}_{i'} \right\rangle \right \vert^{\frac{1}{2}} \right\vert
    \lesssim \frac{1}{\left( 1+ t_1  ( i, i') + t_2 (i, i') \right)^{\frac{\kappa_0}{4}}}
    \cdot \frac{1}{\left(1 + r (i,i') + \theta (i,i') \right)^{\frac{\kappa}{4}}}   
\end{split}
\hspace{0.25cm},
\label{eq:Dritte_Abschaetzung}
\end{equation} which again implies \eqref{eq:Produkt} as long as $\kappa_0$ and $\kappa$ are not smaller than 24. $\Box$
\newline

Theorem 1 indicates that the elements of the Gram matrix $G :=(\langle \psi_i \vert \psi_{i'} \rangle)_                                               {(i,i') \in I \times I}$ of the wave packet system $\{\psi_i\}_{i \in I}$ decay polynomially as they deviate from the main diagonal of the matrix or, in other words, that  the correlation of the wave packets $\psi_i$ and $\psi_{i'}$ fades as their distance $\rho (i, i') $ in the phase space increases. This can be viewed as near orthogonality, so to speak, of the wave packets as their true orthogonality would amount to their Gram matrix being diagonal. The order of the polynome in \eqref{eq:Produkt} can be increased by choosing larger $\kappa_0$ and $\kappa$ and $\kappa_2$. To formulate the next theorem, which is a corollary of Theorem 1 and Lemma 2, we shall need the following definition \cite{Groechenig_2003}.

\begin{definition}
Let $I \subset \mathbb{R}^d$ be a separated denumerable metric index space with a distance $\rho$ and $n > d$. The set of functions $\{ \psi_i \}_{i \in I}$ is said to be $intrinsically$ $localised$ if
\begin{equation}
\vert \langle \psi_i \vert \psi_j \rangle \vert \lesssim (1+ \rho (i, i'))^{-n}
\label{eq:Selbst-Localisation}
\end{equation} for all $i$ and $i' \in I$.
\end{definition}

\noindent \textbf{Theorem 2.} Let $j^{\max}$, $I$ and $W (\alpha, \beta) = \{\psi_i\}_{i \in I}$ be a finite natural number, the metric index space defined by \eqref{eq:Index} with the distance $\rho$ defined by \eqref{eq:rho} and the wave packet system respectively. Then the wave packet system $W (\alpha, \beta)$ is intrinsically localised.
\newline

\noindent $Proof$. The set $I$ is denumerable, according to \eqref{eq:Index}. Once equipped with the distance $\rho$, it becomes a metric space, which will be separated as long as $j^{\max}$ is finite, according to Lemma 2. Furthermore we know from Theorem 1 that the elements of the Gram matrix $G$ of the wave packet system $W (\alpha, \beta) = \{\psi_i\}_{i \in I}$ satisfy \eqref{eq:Selbst-Localisation} with $6 = n > d = 5$ as $i$ and $i' \in I  \subset \mathbb{Z}^5 \subset \mathbb{R}^5$. $\Box$
\newline

Finally from Proposition 1 and 3 in \cite{Jaffard_1990} and Lemma 3 and Theorem 1 of this report we infer that the Gram matrix $G$ of the wave packet system $W (\alpha, \beta) = \{\psi_i\}_{i \in I}$ with $I$ defined by \eqref{eq:Index} with a finite $j^{\max}$ will belong to a solid involutive Banach algebra closed under inversion and therefore the wave packet system will be $self$-$localised$ according to Definition 3 in \cite{Fornasier_2005} of this notion.

\section{Conclusions}

\noindent In this report we first reformulated more neatly the notion of the wave packet system $W (\alpha, \beta)$, originally introduced in \cite{Bytchenkoff_2019}, so that it now forms automatically a Banach frame and a set of atoms of the corresponding wave packet space $\mathcal{W}_{s}^{p, q} (\alpha, \beta)$. We then introduced the notion of the distance between the indexes of the elements of the wave packet system and demonstrated that the index set with this distance is separated. Finally we proved that the absolute value of the scalar product of two elements of the system decays polynomially as the distance between their indexes increases. This, in its turn, can be interpreted as near orthogonality or, in other words, good localisation of the elements the wave packet system in the phase space. We believe that this property will make the wave packet systems particularly useful for the efficient approximation of functions of a wide range of classes and discretisation of bounded linear operators. More specifically, good localisation of the wave packet systems with a wide range of values of the parameters $\alpha$ and $\beta$, namely those that satisfy $0 \leqslant \beta \leqslant \alpha \leqslant 1$ -- which we have now established -- gives one great freedom to choose the system with the type of localisation that would ensure the maximal sparsity of the matrix representing any given Fourier integral operator. In doing so one can be expected to be in a position to further relax the restrictions on the symbol and phase of the Fourier integral operators imposed in the pioneering works \cite{Cordero_2012, Cordero_2013, Cordero_2014} without undermining the sparsity of their representations and to formulate propositions about the boundedness and invertibility of those operators \cite{Borup_2006, Cordero_2013}. 

\section{Notations}

\noindent $\emptyset$ denotes an $empty$ $set$.

\noindent $\bar{A}$ denotes the $adherence$ or $closure$ of the set $A$.

\noindent $A \subset B$ indicates that the set $A$ is a $subset$ of the set $B$.

\noindent $A \cup B$ denotes the $union$ of the sets $A$ and $B$.

\noindent $A \cap B$ denotes the $intersection$ of the sets $A$ and $B$.
\newline

\noindent $B_r(x)$ denotes the Euclidean ball of radius $r$ with its centre in the point $x$.
\newline

\noindent $\mathbb{N}$ denotes the set of all $natural$ $numbers$.

\noindent $\mathbb{N}_0 := \mathbb{N} \cup \{ 0 \}$.

\noindent $\mathbb{N}_0^d$ denotes the set of all $d$-$tuple$ $of$ numbers from $\mathbb{N}_0$.

\noindent $\mathbb{Z}$ denotes the set of all $integer$ $numbers$.

\noindent $\mathbb{Z}^d$ denotes the set of all $d$-$tuple$ $of$ $integer$ $numbers$.

\noindent $\mathbb{R}$ denotes the set of all $real$ $numbers$.

\noindent $\mathbb{R}_+$ denotes the set of all $non$-$negative$ $real$ $numbers$.

\noindent $\mathbb{R}^d$ denotes the set of all $d$-$tuple$ $of$ $real$ $numbers$.

\noindent $\mathbb{C}$ denotes the set of all $complex$ $numbers$.

\noindent $\mathbb{C}^I$ denotes the set of all sequences of complex numbers indexed by $i \in I$.
\newline

\noindent $a \lesssim b$ indicates that there exists such a constant $c$ that $a \leqslant c \, b$.

\noindent $a \asymp b$ indicates that $a \lesssim b$ and $b \lesssim a$.

\noindent $\lceil a \rceil$ denotes the smallest integer number greater than or equal to the real number $a$.
\newline

\noindent The set of all $continuously$ $differentiable$ $functions$ is denoted by $C^1$.

\noindent The set of all $infinitely$ $differentiable$ $functions$ is denoted by $C^\infty$.

\noindent The set of all $infinitely$ $differentiable$ $functions$ $of$ $compact$ $support$ is denoted by $C_c^\infty$.

\noindent The set of all $measurable$ $p$-$integrable$ $functions$ with $p \geqslant 1$ is denoted by $L^p$.

\noindent The set of all $essentially$ $bounded$ $measurable$ $functions$ is denoted by $L^\infty$. 
\newline

\noindent The $composition$ $of$ $functions$ $f$ and $g$ is denoted by
$f \circ g$.
\newline

\noindent The $Fourier$ $transform$ $\hat{f} (\xi)$ of the function $f(t) \in L^1 (\mathbb{R}^d)$ is defined by
\begin{equation}
\hat{f} (\xi) := \left( \mathcal{F} f \right) (\xi) := \int\limits_{\mathbb{R}^d} f(t) \, e^{-2 \pi i \langle t, \xi \rangle} \, dt \hspace{0.5cm} \text{for} \hspace{0.5cm} \xi \in \mathbb{R}^d
\end{equation} where $\langle t, \xi \rangle := \sum_{i = 1}^d t_i \xi_i$. 
\newline

\noindent The $inverse$ $Fourier$ $transform$ $f(t) := \left( \mathcal{F}^{-1} \hat{f} \right) (t)$ of the function $\hat{f} (\xi) \in L^1 (\mathbb{R}^d)$ is defined by
\begin{equation}
f(t) := \left( \mathcal{F}^{-1} \hat{f} \right) (t) := \int\limits_{\mathbb{R}^d} \hat{f} (\xi) \, e^{2 \pi i \langle t, \xi \rangle} \, d \xi \hspace{0.5cm} \text{for} \hspace{0.5cm} t \in \mathbb{R}^d
\hspace{0.25cm}.
\end{equation}
\newline

\noindent $x^a := x_1^{a_1} \cdots x_d^{a_d}$ where $x \in \mathbb{R}_0^d$ and $a \in \mathbb{N}_0^d$.

\noindent $\partial^a f :=
\frac{\partial^{a_d}}{\partial x_d^{a_d}}
\cdots \frac{\partial^{a_1}}{\partial x_1^{a_1}} f$ where $f$ is a differentiable function of $x \in \mathbb{R}_0^d$ and $a \in \mathbb{N}_0^d$.
\newline

\noindent The $Schwartz$ $space$ $\mathcal{S}$ of functions is defined by
\begin{equation}
\mathcal{S} (\mathbb{R}^d)
  := \left\{
        f \in C^\infty (\mathbb{R}^d) \hspace{0.25cm} 
        : \hspace{0.25cm} \vert \vert f \vert \vert_{a, b} < \infty \hspace{0.25cm} \text{for} \hspace{0.25cm} \forall \, a, b \in \mathbb{N}^d_0
     \right\} \, ,
\end{equation} where
\begin{equation}
\vert \vert f \vert \vert_{a, b} :=
                         \sup_{x \in \mathbb{R}^d} \,
                           |x^a ( \partial^b f ) (x)|
                           \hspace{0.25cm}.
\end{equation}

\section*{Acknowledgements}

\noindent I am very grateful to Peter Balazs, Nicki Holighaus, Georg Tauböck and Luis Daniel Abreu for an interesting and helpful discussion we had after I had reported on this work at a seminar at the Acoustics Research Institute of the Austrian Academy of Sciences. I am also grateful to Hans Feichtinger and Felix Voigtlaender for drawing my attention to reports relevant to this work. Finally I thank the Centre National de la Recherche Scientifique and Österreichische Akademie der Wissenschaften for their support of this work.

\end{document}